\documentclass[11pt]{amsart}
\usepackage{amsfonts}
\usepackage{graphicx}
\usepackage{amssymb}
\usepackage{amsmath}
\usepackage{amscd}
\usepackage{mathrsfs}
\usepackage{stmaryrd}
\usepackage{longtable}
\usepackage{txfonts}
\usepackage{wrapfig}
\usepackage{picins}
\usepackage{float}
\usepackage{xypic}
\usepackage{dsfont}
\usepackage{xcolor}
\usepackage{color}
\usepackage[colorlinks,linkcolor=blue,anchorcolor=blue,
citecolor=blue]{hyperref}
\usepackage{hyperref}
\allowdisplaybreaks

\newtheorem{theorem}{Theorem}[subsection]

\newtheorem{lem}[theorem]{Lemma}

\newtheorem{rem}[theorem]{Remark}

\newtheorem{thm}[theorem]{Theorem}

\begin{document}
\setlength{\oddsidemargin}{0cm}
\setlength{\evensidemargin}{0cm}

\title{Classification of three dimensional complex $\omega$-Lie algebras}

\author{Yin Chen}

\address{School of Mathematics and Statistics, Northeast Normal University,
 Changchun 130024, P.R. China}

\email{ychen@nenu.edu.cn}

\author{Chang Liu}

\address{School of Mathematics and Statistics, Northeast Normal University,
 Changchun 130024, P.R. China}

\email{liuc813@nenu.edu.cn}

\author{Run-Xuan Zhang}

\address{School of Mathematics and Statistics, Northeast Normal University,
 Changchun 130024, P.R. China}

\email{zhangrx728@nenu.edu.cn}

\date{\today}

\def\shorttitle{Classification of three dimensional complex $\omega$-Lie algebras}

\begin{abstract}
A complex $\omega$-Lie algebra is a  vector space $L$ over the complex field, equipped with a skew symmetric bracket $[-,-]$ and a bilinear form
$\omega$ such that $$[[x,y],z]+[[y,z],x]+ [[z,x],y]=\omega(x,y)z+\omega(y,z)x+\omega(z,x)y$$ for all $x,y,z\in L$. The notion of $\omega$-Lie algebras, as a generalization of Lie algebras, was introduced in Nurowski \cite{Nur2007}. Fundamental results about finite-dimensional  $\omega$-Lie algebras were developed by Zusmanovich \cite{Zus2010}. In \cite{Nur2007}, all three dimensional non-Lie real $\omega$-Lie algebras were classified. The purpose of this note is to provide an approach to classify all three dimensional non-Lie complex $\omega$-Lie algebras. Our method also gives a new proof of the classification in Nurowski \cite{Nur2007}.
\end{abstract}



\maketitle
\baselineskip=19pt

\setcounter{subsection}{0}
\renewcommand{\thesubsection}
{\arabic{subsection}}

\setcounter{equation}{0}
\renewcommand{\theequation}
{\arabic{equation}}

\setcounter{theorem}{0}
\renewcommand{\thetheorem}
{\arabic{theorem}}

\subsection{Introduction}

Let $k$ be a field of characteristic zero and $L$ be a finite-dimensional vector space over $k$. Let $[-.-]:L\times L\longrightarrow L$ be an anti-commutative product on $L$ and $\omega:L\times L\longrightarrow k$ be a bilinear form on $L$.
The triple $(L,[-.-],\omega)$ is called an $\omega$-\textit{Lie algebra} if the following condition is satisfied:
\begin{equation}\label{Jacob-identity}
[[x,y],z]+ [[y,z],x]+  [[z,x],y]=\omega(x,y)z+\omega(y,z)x+\omega(z,x)y
\end{equation}
for all $x,y,z\in L$. The equation (\ref{Jacob-identity}) is called the $\omega$-\textit{Jacobi identity}.
Apparently, the $\omega$ is also skew-symmetric; an $\omega$-Lie algebra is a  Lie algebra if and only if the bilinear form $\omega\equiv 0$.
So we usually call the Lie algebras \textit{trivial} $\omega$-Lie algebras.

The notion of $\omega$-Lie algebras, which is related to the study of isoparametric hypersufaces in Riemannian geometry (\cite{BN2007,Nur2008}), was introduced in the recent work of Nurowski \cite{Nur2007}. By the definition, it is easy to see that all $\omega$-Lie algebras are trivial in the case of dimension 1 and 2.
The first example of  nontrivial 3-dimensional $\omega$-Lie algebra was given in \cite{Nur2007}.
In that paper, Nurowski finally  completed the classification of 3-dimensional $\omega$-Lie algebra over the field of real numbers.

A fundamental development of  $\omega$-Lie algebras was by Zusmanovich \cite{Zus2010}, in which a lot of basic concepts, such as modules, (quasi-) ideals and (generalized) derivations, were introduced; some basic properties of $\omega$-Lie algebras were found.
One of Zusmanovich's results asserts that finite-dimensional nontrivial  $\omega$-Lie algebras are either low-dimensional or have
an abelian  subalgebra of small codimension with some restrictive conditions. In particular, the following useful result is proved.

\begin{lem}[\cite{Zus2010}, Lemma 8.1]\label{Zus2010}
If $L$ is a finite-dimensional $\omega$-Lie algebra with non-degenerate $\omega$, then dim $L=2$.
\end{lem}

Recall that a skew-symmetric bilinear form $\omega$ on $L$ is \textit{degenerate} if there exists a nonzero vector $x\in L$ such that $\omega(x,y)=0$ for all $y\in L$. Lemma \ref{Zus2010} mentioned above means that the bilinear form $\omega$ on any $\omega$-Lie algebra $L$ must be degenerate if dim$L\geq 3$.

The main purpose of this note is to provide an approach to classify all 3-dimensional  nontrivial $\omega$-Lie algebras over the field of complex numbers. In final remarks (Section 5), our method is also applied to give a classification of all 3-dimensional  nontrivial real $\omega$-Lie algebras.

In what follows, $\mathbb{R}$ and $\mathbb{C}$ are the fields of real and complex numbers respectively, and $L$ denotes an $\omega$-Lie algebra with a basis $\{x,y,z\}$.
We write $\wedge^{2} L$ for the exterior power of $L$ with the basis $\{x\wedge y,x\wedge z,y\wedge z\}$ and $\varphi=[-,-]:\wedge^{2} L\longrightarrow L$ is the bracket product. We use $L'=[L,L]$ to denote the commutator subalgebra of $L$.
We call the dimension of $L'$ the \textit{rank} of $\varphi$. In the book \cite{FH1991},
Fulton and  Harris  presented a classification of 3-dimensional  complex Lie algebras by considering the rank of $\varphi$ case by case.

In this note, we will follow some ideas in \cite{FH1991} and discuss the rank of $\varphi$, which may be $0,1,2,3$.

The following theorem is our main result.

\begin{thm}\label{main-thm}
Let $L$ be a nontrivial (i.e. non-Lie) 3-dimensional $\omega$-Lie algebra, then it must be isomorphic to one of the following algebras:
\begin{enumerate}
  \item \quad $L_{1}:\quad [x,z]=0,[y,z]=z, [x,y]=y\textrm{ and }\omega(y,z)=\omega(x,z)=0,\omega(x,y)=1.$
  \item \quad $L_{2}:\quad [x,y]=0,[x,z]=y,[y,z]=z \textrm{ and } \omega(x,y)=0, \omega(x,z)=1, \omega(y,z)=0.$
  \item \quad
  $A_{\alpha} :\quad [x,y]=x,[x,z]=x+y, [y,z]=z+\alpha x\textrm{ and  }
\omega(x,y)=\omega(x,z)=0,$
\begin{center}
$\omega(y,z)=-1,\textrm{ where }\alpha\in \mathbb{C}.$
\end{center}
  \item \quad $B: \quad [x,y]=y, [x,z]=y+z, [y,z]=x\textrm{ and } \omega(x,y)=\omega(x,z)=0,$
$\omega(y,z)=2.$
   \item \quad $C_{\alpha} : \quad [x,y]=y, [x,z]=\alpha z, [y,z]=x\textrm{ and } \omega(x,y)=\omega(x,z)=0,$
  \begin{center}
$\omega(y,z)=1+\alpha,\textrm{ where }0\neq\alpha\in \mathbb{C}.$
  \end{center}
\end{enumerate}
\end{thm}

\subsection{Ranks 0 and 1}
We continue to follow the notations in the preceding section.

If the rank of $\varphi$ is zero, then $L$ is abelian. It follows from the $\omega$-Jacobi identity (\ref{Jacob-identity}) that
$$\omega(x,y)z+\omega(y,z)x+\omega(z,x)y=0.$$
Since $\{x,y,z\}$ is a basis of $L$, $\omega(x,y)=\omega(y,z)=\omega(z,x)=0.$ Thus in this case, the $\omega$-Lie structure on $L$ is trivial.

If the rank of $\varphi$ is 1,
then dim$L'=1$ and
the kernel of  $\varphi$ is two dimensional. Suppose that $\{x,y,z\}$ is a basis of $L$ such that $[x,y]=[x,z]=0$. We let $[y,z]=ax+by+cz$ for some $a,b,c\in \mathbb{C}$.
By $\omega$-Jacobi identity, we have
\begin{eqnarray*}
&& \omega(x,y)z+\omega(y,z)x+\omega(z,x)y\\
 &=&[[x,y],z]+ [[y,z],x]+  [[z,x],y]\\
 &=& [ax+by+cz,x]=0.
\end{eqnarray*}
Since $x,y,z$ is linearly independent,  $\omega(x,y)=\omega(y,z)=\omega(z,x)=0$. Thus  $\omega$ is trivial.
The same arguments as in \cite{FH1991} (page 137) will imply that there only exist two Lie algebras:
\begin{eqnarray*}
\mathfrak{g}_{1} &:& [x,y]=[x,z]=0,\textrm{ and }[y,z]=y. \\
\mathfrak{g}_{2}  &:& [x,y]=[x,z]=0,\textrm{ and } [y,z]=x.
\end{eqnarray*}
They are trivial $\omega$-Lie algebras.

\subsection{Rank 2}
In this case,
we choose $\{y,z\}$ as a basis of $L'$ and $x\notin L'$. We assume that $[y,z]=ay+bz$. Our arguments will be separated into the following two cases:
$a=b=0$ or the others.

\textcolor[rgb]{0.00,0.07,1.00}{Case 1. }\quad If both $a$ and $b$ are zero, then $[x,y]\neq 0$ and
$[x,z]\neq 0$ because the kernel of  $\varphi$ is one dimensional.
Thus the linear map ad$_{x}:L'\longrightarrow L'$ by $u\longmapsto [x,u]$ is bijective. By linear algebra, we can choose the suitable basis elements $y,z$ such that ad$_{x}$
is similar to
$$\left(
  \begin{array}{cc}
    c & 0 \\
    0 & d \\
  \end{array}
\right)~~\textrm{or}~~\left(
  \begin{array}{cc}
    e & 0 \\
    1 & e \\
  \end{array}
\right),$$
where $c,d,e$ are nonzero complex numbers. In the first situation, $[x,y]=\textrm{ad}_{x}(y)=cy$ and $[x,z]=\textrm{ad}_{x}(z)=dz$.
Let $\widetilde{x}=c^{-1}x$, then $[\widetilde{x},y]=y$ and $[\widetilde{x},z]=(c^{-1}d)z$. We obtain a family of Lie algebras with one parameter:
\begin{eqnarray*}
\mathfrak{h}_{\alpha} &:& [x,y]=y,[x,z]=\alpha z\textrm{ and } [y,z]=0,\textrm{ where }0\neq\alpha\in \mathbb{C}.
\end{eqnarray*}

In the second situation, $[x,y]=ey$ and $[x,z]=y+ez$. Let Let $\widetilde{x}=e^{-1}x$, then $[\widetilde{x},y]=y$ and
$[\widetilde{x},z]=e^{-1}y+z$. Let $\widetilde{y}=e^{-1}y$, then $[\widetilde{x},\widetilde{y}]=\widetilde{y}$ and $[\widetilde{x},z]=\widetilde{y}+z$. Thus we get a Lie algebra:
\begin{eqnarray*}
\mathfrak{g}_{3} &:& [x,y]=y,[x,z]=y+z\textrm{ and }[y,z]=0.
\end{eqnarray*}

\textcolor[rgb]{0.00,0.07,1.00}{Case 2.} \quad Assume that one of  $a$ and $b$ is not zero.  We need only to consider the case of $b\neq 0$ because if $a\neq 0$ then we can transpose
$y$ and $z$ and will get the same results.

Let $\widetilde{z}=z+b^{-1}ay$, then $[y,\widetilde{z}]=[y,z]=ay+bz=b\widetilde{z}$. Let $\widetilde{y}=b^{-1}y$ then $[\widetilde{y},\widetilde{z}]=\widetilde{z}$.
 So in this case, we can assume that $[y,z]=z$.
Since the  kernel of  $\varphi$ is one dimensional, one of $[x,z]$ and $[x,y]$ is zero, and the other is not zero.

Subcase 1. \quad If $[x,y]=ay+cz\neq 0$ and $[x,z]=0$, then
\begin{eqnarray*}
&& \omega(x,y)z+\omega(y,z)x+\omega(z,x)y\\
 &=&[[x,y],z]+ [[y,z],x]+  [[z,x],y]\\
 &=& [ay+cz,z]+[z,x]=az.
\end{eqnarray*}
This means that $\omega(y,z)=\omega(z,x)=0$ and $a=\omega(x,y)$. If $c=0$, then $a\neq0$. We let $\widetilde{x}=a^{-1}x$ and it is easy to check that
$[\widetilde{x},z]=0,[y,z]=z$, $[\widetilde{x},y]=y$, $\omega(y,z)=\omega(\widetilde{x},z)=0$ and $\omega(\widetilde{x},y)=1$.
Thus we get a nontrivial $\omega$-Lie algebra:
$$
L_{1}:\quad [x,z]=0,[y,z]=z, [x,y]=y\textrm{ and  }\omega(y,z)=\omega(x,z)=0,\omega(x,y)=1.
$$
If $c\neq 0$ and $a=0$, then the $\omega$-Jacobi identity implies that $\omega$ is trivial. If $c\neq 0$ and $a\neq 0$, we assume that
$\widetilde{y}=y+a^{-1}cz$, then $[x,z]=0,[\widetilde{y},z]=z,[x,\widetilde{y}]=a\widetilde{y}$ and $\omega(\widetilde{y},z)=\omega(x,z)=0,\omega(x,\widetilde{y})=a$, so it is easy to see that the corresponding  $\omega$-Lie algebra  is isomorphic to $L_{1}$.

Subcase 2. \quad If $[z,x]=-[x,z]=ay+cz\neq 0$ and $[x,y]=0$, then
\begin{eqnarray*}
&& \omega(x,y)z+\omega(y,z)x+\omega(z,x)y\\
 &=&[[x,y],z]+ [[y,z],x]+ [[z,x],y]\\
 &=& [z,x]+[cz,y]=ay.
\end{eqnarray*}
Thus $\omega(y,z)=0$, $\omega(x,y)=0$ and $\omega(z,x)=a$. Now we have
$$[y,z]=z,[x,y]=0,[z,x]=ay+cz.$$
Since the dimension of $L'$ is two, $a\neq 0$.
Let $\widetilde{y}=y+a^{-1}cz$, then
 $$[z,x]=a\widetilde{y}, [\widetilde{y},z]=z,[x,\widetilde{y}]=-c\widetilde{y}.$$
Recall that the  kernel of  $\varphi$ is one dimensional, so $c$ must be zero.
Now we can assume that $\widetilde{z}=a^{-1}z$ and $\widetilde{x}=-x$. This yields that
$$[\widetilde{y},\widetilde{z}]=\widetilde{z},[\widetilde{x},\widetilde{y}]=0,[\widetilde{x},\widetilde{z}]=\widetilde{y}; \omega(\widetilde{y},\widetilde{z})=0, \omega(\widetilde{x},\widetilde{y})=0,\omega(\widetilde{z},\widetilde{x})=-1.$$
Hence there is a nontrivial $\omega$-Lie algebra:
$$L_{2}:\quad[x,y]=0,[x,z]=y,[y,z]=z\textrm{ and }\omega(x,y)=0, \omega(x,z)=1,\omega(y,z)=0.$$

\subsection{Rank 3}

For the case of $\omega\equiv 0$, it follows from the Fulton and  Harris's arguments in \cite{FH1991} (pages 141-142) that there exists only one Lie algebra:
\begin{eqnarray*}
\mathfrak{g}_{4} &:& [x,y]=2y,[x,z]=-2z\textrm{ and }[y,z]=x.
\end{eqnarray*}

Next we consider the nontrivial case.
Since the dimension of $L'$ is 3,  the rank of adjoint map ad$_{x}:L\longrightarrow L$ must be 2 for any nonzero $x\in L$.
Thus the kernel of ad$_{x}$ is equal to $\mathbb{C}\cdot x$.

By Lemma \ref{Zus2010}, if $\omega$ is non-degenerate, then $L$ must have dimension 2. So  the bilinear form $\omega$ we consider here is degenerate. This means that there exists an nonzero element $x\in L$ such that $\omega(x,v)=0$ for all $v\in L$.
Now we fix $x$. By the Jordan canonical form, we can choose a suitable basis $\{u,y,z\}$ of $L$ such that ad$_{x}$ is similar to
\begin{equation}\label{jordan}
A=\left(
    \begin{array}{ccc}
      0& 0 & 0 \\
      1 & 0 & 0 \\
      0 & 1 & 0 \\
    \end{array}
  \right),\quad B=\left(
    \begin{array}{ccc}
      0& 0 & 0 \\
      0 & \delta & 0 \\
      0 & 1 & \delta \\
    \end{array}
  \right),\quad C=\left(
    \begin{array}{ccc}
      0& 0 & 0 \\
      0 & \mu & 0 \\
      0 & 0 & \nu \\
    \end{array} \right),\textrm{ or } D=\left(
    \begin{array}{ccc}
      0& 0 & 0 \\
      1 & 0 & 0 \\
      0 & 0 & \tau \\
    \end{array} \right),
\end{equation}
where $\delta,\mu,\nu$ and $\tau$ are all nonzero and are the eigenvalues of ad$_{x}$.
Thus our arguments consist of the following four cases.

\textcolor[rgb]{0.00,0.07,1.00}{Case 1. }\quad If  ad$_{x}$ is similar to $A$, then $[x,u]=0,[x,y]=u$ and $[x,z]=y$. Let $u=ax$ for some nonzero $a\in \mathbb{C}$.
Then we have $[x,y]=ax$ and $[x,z]=y.$ That is, $[a^{-1}x,a^{-1}y]=a^{-1}x$ and $[a^{-1}x,z]=a^{-1}y$.
Thus we can assume that $\{x,y,z\}$ is a basis of $L$ such that
$$[x,y]=x, [x,z]=y.$$
Next we need to determine the commutator relations of $y$ and $z$.
Let $[y,z]=bx+cy+dz$, then
\begin{eqnarray*}
&& \omega(x,y)z+\omega(y,z)x+\omega(z,x)y\\
 &=&[[x,y],z]+ [[y,z],x]+ [[z,x],y]\\
 &=& [x,z]+[cy+dz,x]\\
 &=& (1-d)y-cx.
\end{eqnarray*}
Thus $\omega(x,y)=0, \omega(x,z)=d-1, \omega(y,z)=-c$. Since $\omega(x,z)=0$, $d=1$. Notice that $\omega$ is not trivial, so $c\neq 0$.
In the equation $[y,z]=bx+cy+z$,
we first can assume that $\widetilde{z}=z+cy$. Then
$$[y,\widetilde{z}]=\widetilde{z}+bx,[x,y]=x,[x,\widetilde{z}]=cx+y,$$
and $\omega(x,y)=0, \omega(x,\widetilde{z})=0$ and $\omega(y,\widetilde{z})=-c$.
Let $\alpha=c^{-1}b,\beta=c^{-1},z'=c^{-1}\widetilde{z}$, then
$$[y,z']=z'+\alpha x,[x,y]=x,[x,z']=x+\beta y$$
and $\omega(x,y)=0, \omega(x,z')=0$ and $\omega(y,z')=-1$.
We define $x'=\beta^{-1}x$ and $\gamma=\alpha\beta$. Then
$$[x',y]=x',[x',z']=x'+y, [y,z']=z'+\gamma x'$$
with $ \omega(x',y)=0, \omega(x',z')=0, \omega(y,z')=-1$.
Hence we get a family of $\omega$-Lie algebras with one parameter:
\begin{eqnarray*}
A_{\alpha} &:& [x,y]=x,[x,z]=x+y, [y,z]=z+\alpha x;\\
&& \omega(x,y)=0, \omega(x,z)=0, \omega(y,z)=-1,\textrm{ where }\alpha\in \mathbb{C}.
\end{eqnarray*}

\textcolor[rgb]{0.00,0.07,1.00}{Case 2. }\quad If  ad$_{x}$ is similar to $B$, then $[x,u]=0, [x,y]=\delta y,[x,z]=y+\delta z$.
Obviously, $[\delta^{-1} x,\delta^{-1} y]=\delta^{-1} y,[\delta^{-1} x,z]=\delta^{-1} y+ z$.
Thus we can assume that $\{x,y,z\}$ is a basis of $L$ such that
$$[x,y]=y, [x,z]=y+z.$$
Let $[y,z]=ax+by+cz$, then
\begin{eqnarray*}
&& \omega(x,y)z+\omega(y,z)x+\omega(z,x)y\\
 &=&[[x,y],z]+ [[y,z],x]+ [[z,x],y]\\
 &=& 2[y,z]+[by+cz,x]\\
 &=& 2a x+(b-c)y+c z.
\end{eqnarray*}
Thus $\omega(x,y)=c, \omega(x,z)=c-b$ and $\omega(y,z)=2a$. Recall that $x$ belongs to the kernel of $\omega$, so $c=b=0$.
Hence we have
\begin{eqnarray*}
&& [x,y]=y, [x,z]=y+z, [y,z]=\alpha x;\\
&& \omega(x,y)=0, \omega(x,z)=0, \omega(y,z)=2\alpha,\textrm{ where }0\neq\alpha\in \mathbb{C}.
\end{eqnarray*}
Let $u=(\sqrt{\alpha})^{-1}$, then $$[x,uy]=uy, [x,uz]=uy+uz, [uy,uz]=x.$$
Therefore, we obtain a new $\omega$-Lie algebra:
\begin{eqnarray*}
B&:& [x,y]=y, [x,z]=y+z, [y,z]=x;\\
&& \omega(x,y)=0, \omega(x,z)=0, \omega(y,z)=2.
\end{eqnarray*}

\textcolor[rgb]{0.00,0.07,1.00}{Case 3. }\quad If  ad$_{x}$ is similar to $C$, then
$[x,u]=0, [x,y]=\mu y,[x,z]=\nu  z$. We can assume that $\{x,y,z\}$ is a basis of $L$. Let $[y,z]=ax+by+cz$.
Notice that the dimension of $L'$ is 3, so $a\neq 0$. Thus as before, we can assume that
$[x,y]=y,[x,z]=\alpha  z$ and $[y,z]=x+by+cz$, where $\alpha\neq 0$.
The $\omega$-Jacobi identity implies that
$$b=c=0\textrm{ and }\omega(x,y)=0, \omega(x,z)=0, \omega(y,z)=1+\alpha.$$
We obtain again  a family of $\omega$-Lie algebras of one parameter:
\begin{eqnarray*}
C_{\alpha} &:& [x,y]=y, [x,z]=\alpha  z, [y,z]=x;\\
&& \omega(x,y)=0, \omega(x,z)=0, \omega(y,z)=1+\alpha,\textrm{ where }0\neq\alpha\in \mathbb{C}.
\end{eqnarray*}

\textcolor[rgb]{0.00,0.07,1.00}{Case 4. }\quad If  ad$_{x}$ is similar to $D$, then
$[x,u]=0, [x,y]=u,[x,z]=\tau  z$. Let $\{x,y,z\}$ be a basis of $L$ and $[y,z]=ax+by+cz$.
As before, we can assume that
$[x,y]=x,[x,z]=z$ and $[y,z]=ax+by+cz$.
By the $\omega$-Jacobi identity we will see that $\omega(x,y)=1$, which contradicts with the assumption that $\omega(x,v)=0$ for all $v\in L$. Thus in this situation, there do not exist any $\omega$-Lie algebras.

This completes the proof of Theorem \ref{main-thm}.

\subsection{Remarks}

We close this paper with the following remarks, in which almost all of the ideas come from the referee's report on the first version of our paper.

The method that we use in this paper can actually be  applied to classify all three dimensional real $\omega$-Lie algebras.
The discussion of  rank$(\varphi)=0$ or 1 in Section 2 can be carried to the case of real numbers, and we do not obtain any new $\omega$-Lie algebras.

Notice that our proof of Theorem \ref{main-thm} relies deeply on the Jordan canonical form (see Sections 3 and 4) and the fact that each nonzero complex number has a square root in $\mathbb{C}$. Moreover,
we are interested in nontrivial (i.e., non-Lie) $\omega$-Lie algebras. So in the case of rank$(\varphi)=2$, the direct computation shows that there are not any new nontrivial $\omega$-Lie algebras. However, in the case of rank 3,  some interesting things will happen. First, two families of new nontrivial $\omega$-Lie algebras (not contained in the list of Theorem \ref{main-thm}) will appear, because there are five possibilities of the similar canonical form. Second, the canonical form $B$ will give two real non-isomorphic $\omega$-Lie algebras.

In Remark \ref{rem1} below, we explain the reason why the complex $\omega$-Lie algebra $B$ in Theorem \ref{main-thm} can ``split" into two real $\omega$-Lie algebras $\mathcal{B}_{1}$ and $\mathcal{B}_{-1}$.
 We also describe $\mathcal{E}_{\alpha}^{+}$ and $\mathcal{E}_{\alpha}^{-}$, the two new families of nontrivial $\omega$-Lie algebras over $\mathbb{R}$ in details. Thus we write out in advance all nontrivial real $\omega$-Lie algebras up to isomorphism:
 \begin{equation}\label{real}
L_{1}, L_{2}, \mathcal{A}_{\alpha}, \mathcal{C}_{\alpha}, \mathcal{B}_{1}, \mathcal{B}_{-1}, \mathcal{E}_{\alpha}^{+},\mathcal{E}_{\alpha}^{-}.
 \end{equation}
Another reason why Nurowski's list is bigger than the above (\ref{real}) is that  Nurowski considered the action of orthogonal group $O(\mathbb{R})$ and classified all $\omega$-Lie algebras up to ``equivalence". In Remark \ref{rem2}, we compare Nurowski's list with ours.

\begin{rem}\label{rem1}
{\rm
We employ the notations in Section 4, and assume that the ground field in the consideration is $\mathbb{R}$. By linear algebra, we can choose a suitable basis $\{u,y,z\}$ of $L$ such that
$$\omega(x,v)=0,\forall v\in L,$$ and
ad$_{x}$ is similar to the canonical forms $A,B,C,D$ in (\ref{jordan}), or
$$E=\left(
    \begin{array}{ccc}
      0& 0 & 0 \\
      0 & 0 & -b \\
      0 & 1 & a \\
    \end{array}
  \right),
$$ where $a^{2}<4b$.

In the cases $A,C,D$, we proceed as in the case over $\mathbb{C}$ before. Thus we will obtain two families of real $\omega$-Lie algebras: $\mathcal{A}_{\alpha}$ and $\mathcal{C}_{\alpha}$, which correspond to $A_{\alpha}$ and $C_{\alpha}$ respectively in our Theorem \ref{main-thm}, with the small difference that the parameter $\alpha$ belongs to $\mathbb{R}$ not $\mathbb{C}$.

The case when ad$_{x}$ is similar to $B$ is very interesting because it gives one of the reasons why
our list in Theorem \ref{main-thm} can split into Nurowski's list. First of all, as in Case 2 of Section 4, we obtain a family of $\omega$-Lie algebras
\begin{eqnarray*}
B_{\alpha}&:& [x,y]=y, [x,z]=y+z, [y,z]=\alpha x;\\
&& \omega(x,y)=0, \omega(x,z)=0, \omega(y,z)=2\alpha,\textrm{ where }0\neq\alpha\in \mathbb{R}.
\end{eqnarray*}
There are two cases: $\alpha>0$ and $\alpha<0$. With similar arguments, we obtain two new nontrivial  $\omega$-Lie algebras:
\begin{eqnarray*}
\mathcal{B}_{1}&:& [x,y]=y, [x,z]=y+z, [y,z]=x;\\
&& \omega(x,y)=0, \omega(x,z)=0, \omega(y,z)=2,\quad (\textrm{when } \alpha>0.)
\end{eqnarray*}
and
\begin{eqnarray*}
\mathcal{B}_{-1}&:& [x,y]=y, [x,z]=y+z, [y,z]=-x;\\
&& \omega(x,y)=0, \omega(x,z)=0, \omega(y,z)=-2.\quad (\textrm{when } \alpha<0.)
\end{eqnarray*}
It is not difficult to show that they are isomorphic as complex $\omega$-Lie algebras but not as real $\omega$-Lie algebras.

The case when ad$_{x}$ is similar to $E$ is new, and here we discuss it in detail.
Notice that $b>0$.  We let $u=\sqrt{b}$ and consider the following matrix:
$$E'=\left(
    \begin{array}{ccc}
      0& 0 & 0 \\
      0 & 0 & -u \\
      0 & u & a \\
    \end{array}
  \right),
$$ which is similar to $E$ because they have the same characteristic polynomial. We choose a suitable basis $\{x,y,z\}$ of $L$ such that $\omega(x,L)=0$, and
 ad$_{x}$ is similar to $E'$. Thus
\begin{eqnarray*}
[x,y]&=& -u z, \\
~[x,z]&=& u y+az.
\end{eqnarray*}
So $[u^{-1}x,u y]= -u z$ and $[u^{-1}x,uz]= u y+au^{-1}(uz)$. Therefore we can assume that $\{x,y,z\}$ is a basis of $L$ such that
$$[x,y]=-z, [x,z]=y+\alpha z, $$ where $\alpha\in \mathbb{R}$.
Now we need only to determine the commutator relations of $y$ and $z$.
As before, we let $[y,z]=cx+dy+ez$. Since $\omega(x,y)=0, \omega(x,z)=0$,  the $\omega$-Jacob identity implies that $d=e=0,c\neq 0$ and $\omega(y,z)=\alpha c\neq 0$. Thus we have
\begin{eqnarray*}
&& [x,y]=-z, [x,z]=y+\alpha  z, [y,z]=c x;\\
&& \omega(x,y)=0, \omega(x,z)=0, \omega(y,z)=c\alpha,\textrm{ where }0\neq c,\alpha\in \mathbb{R}.
\end{eqnarray*}
Notice that if $\alpha=0$ then we  obtain a trivial $\omega$-Lie algebra. Thus $\alpha\neq0$.

If $c>0$, we write $v=(\sqrt{c})^{-1}>0$. Thus
$$[x,v y]=-v z, [x,v z]=v y+(\alpha\sqrt{c}) v z, [vy,v z]=x.$$
We replace  $v y$, $v z$, $\alpha\sqrt{c}$ by $y$, $z$ and $\alpha$  respectively, and finally obtain a family of new real $\omega$-Lie algebras with one parameter:
\begin{eqnarray*}
\mathcal{E}_{\alpha}^{+} &:& [x,y]=-z, [x,z]=y+\alpha  z, [y,z]=x;\\
&& \omega(x,y)=0, \omega(x,z)=0, \omega(y,z)=\alpha,\textrm{ where }0\neq\alpha\in \mathbb{R}.
\end{eqnarray*}

If $c<0$, we proceed in a similar way and  obtain another family of $\omega$-Lie algebras:
\begin{eqnarray*}
\mathcal{E}_{\alpha}^{-} &:& [x,y]=-z, [x,z]=y+\alpha  z, [y,z]=-x;\\
&& \omega(x,y)=0, \omega(x,z)=0, \omega(y,z)=-\alpha,\textrm{ where }0\neq\alpha\in \mathbb{R}.
\end{eqnarray*}

}\end{rem}

\begin{rem}\label{rem2}{\rm For the convenience of the reader,
we  reproduce  Nurowski's list in Theorem 2.1 of  \cite{Nur2007} in such a way that it corresponds to our results more clearly.
We consider all $\omega$-Lie algebra of Bianchi type case by case. In order to illustrate the sort of arguments that we use to prove that an $\omega$-Lie algebra of Bianchi type is isomorphic to some  in  our list, we present an example as follows.

We choose  $IV_{T}$.
According to the commutation relations given in the list of Nurowski, we observe that $\omega(x,L)=0$ and ad$x=\begin{pmatrix}
                                                                                                                 0 & 0 &0\\
                                                                                                                 0 & 1 &0\\
                                                                                                                 0 & 0 &1\\
                                                                                                               \end{pmatrix}
$, which is just the matrix $C$ in (\ref{jordan}) with $\mu=\nu=1$. Thus $IV_{T}$ must be isomorphic to some $\omega$-Lie algebra belonging to the family $\mathcal{C}_{\alpha}$. Actually, the direct computation shows that $IV_{T}\simeq \mathcal{C}_{1}$.

Analogous routes and techniques can be applied to other cases, so
we here only list all the correspondence between Nurowski's list and ours (see Table 1 below), and  the detailed proofs are not given.

\begin{table}[!h]
\tabcolsep 5mm
\caption{Nontrivial Real Three-dimensional  $\omega$-Lie Algebras}
  \centering
  \begin{tabular}{c|c|c}
    \hline
Bianchi type & Commutation relations & Being isomorphic to \\ \hline
 $IV_{T}$ & $\begin{array}{c}
    [x,y]=y,[x,z]=z,[y,z]=x\\
    \omega(x,y)=\omega(x,z)=0, \omega(y,z)=-2 \\
  \end{array}$& $\mathcal{C}_{1}$ \\  \hline
 $VI_{T}$ &   $\begin{array}{c}
    [x,y]=y,[x,z]=y+z,[y,z]=x\\
    \omega(x,y)=\omega(x,z)=0, \omega(y,z)=-2 \\
  \end{array}$ & $\mathcal{B}_{1}$ \\  \hline
 $VI_{S}$ &   $\begin{array}{c}
    [x,y]=-x,[x,z]=y,[y,z]=x+z\\
    \omega(x,y)=\omega(y,z)=0, \omega(x,z)=-2 \\
  \end{array}$ & $\mathcal{B}_{1}$ \\  \hline
 $VI_{N}$ &  $\begin{array}{c}
    [x,y]=-x+y,[x,z]=y+z,[y,z]=x+z\\
    \omega(x,y)=0,\omega(y,z)=\omega(x,z)=-2 \\
  \end{array}$ & $L_{1}$ \\  \hline
 $VII_{T}$ &  $\begin{array}{c}
    [x,y]=y,[x,z]=-y+z,[y,z]=x\\
    \omega(x,y)=\omega(x,z)=0,\omega(y,z)=-2 \\
  \end{array}$ & $\mathcal{B}_{-1}$ \\  \hline
  $VIII_{a}$ & $\begin{array}{c}
    [x,y]=-z,[x,z]=-ax-y,[y,z]=x-ay\\
    \omega(x,y)=2a,\omega(x,z)=\omega(y,z)=0, (a>0) \\
  \end{array}$ & $\mathcal{E}_{\alpha}^{+}$ \\  \hline
  $VIII_{T_{a}}$ & $\begin{array}{c}
    [x,y]=ay-z,[x,z]=-y+az,[y,z]=x\\
    \omega(x,y)=\omega(x,z)=0,\omega(y,z)=-2a, (a>0) \\
  \end{array}$ & $\begin{array}{c}
    L_{2}, \textrm{ if }a=1\\
    \mathcal{C}_{\alpha},\textrm{ if }a\neq 1\\
  \end{array}$ \\  \hline
  $VIII_{N_{a}}$ & $\begin{array}{c}
    [x,y]=ay-z,[x,z]=-ax-y+az,[y,z]=x-ay\\
    \omega(x,y)=2a,\omega(x,z)=0,\omega(y,z)=-2a, (a>0) \\
  \end{array}$ & $\mathcal{A}_{\alpha}$ \\  \hline
  $IX_{a}$ & $\begin{array}{c}
    [x,y]=z,[x,z]=-ax-y,[y,z]=x-ay\\
    \omega(x,y)=-2a,\omega(x,z)=\omega(y,z)=0, (a>0) \\
  \end{array}$ & $\mathcal{E}_{\alpha}^{-}$ \\  \hline
  \end{tabular}
\end{table}
}\end{rem}

\subsection*{{\it Acknowledgments}}
The authors are grateful to an anonymous referee for his/her valuable comments
and suggestions on the first version of our paper. In particular, the origin of all presentations in Section 5 comes from the referee's report.
This work was supported partially by NSF of China (11301061, 11226051), STDRP of Jilin Province (20130522098JH, 20140520052JH), and IETPUS (201210200041).

\end{document}